\documentclass[12pt]{article}

\usepackage[english]{babel}
\usepackage[T1]{fontenc}
\usepackage{amsmath}
\usepackage{amsfonts}
\usepackage{amssymb}
\usepackage[utf8]{inputenc}
\usepackage{graphicx,color,url}
\usepackage[margin=2.4cm]{geometry}
\usepackage{booktabs}
\usepackage{cite}
\usepackage{caption}
\usepackage{tikz}
\newtheorem{rem}{Remark}

\def\mca{\mathcal{A}}
\def\mcq{\mathcal{Q}}
\def\wmca{\widetilde{\mathcal{A}}}
\def\wmcq{\widetilde{\mathcal{Q}}}

\def\pmatrix{\left(\begin{array}}
\def\endpmatrix{\end{array}\right)}

\newcommand{\upd}{\,\mathrm{d}}

\author{G.\,Frasca-Caccia\,\quad P.\,E.\,Hydon\,\\[.5cm]
\small
School of Mathematics, Statistics and Actuarial Science\\
\small
University of Kent, Canterbury, CT2 7NZ\\
\small}

\begin{document}

\title{Locally conservative finite difference schemes for the modified KdV equation}

\author{
{\sc G.\,Frasca-Caccia,\,\quad P.\,E.\,Hydon\thanks{Corresponding author. Email: P.E.Hydon@kent.ac.uk}} \\[2pt]
School of Mathematics, Statistics and Actuarial Science\\
University of Kent, Canterbury, CT2 7FS, UK}

\maketitle

\begin{abstract} {

Finite difference schemes that preserve two conservation laws of a given partial differential equation can be found directly by a recently-developed symbolic approach. Until now, this has been used only for equations with quadratic nonlinearity.

In principle, a simplified version of the direct approach also works for equations with polynomial nonlinearity of higher degree. For the modified Korteweg-de Vries equation,
whose nonlinear term is cubic, this approach yields several new families of second-order accurate schemes that preserve mass and either energy or momentum. Two of these families contain Average Vector Field schemes of the type developed by Quispel and co-workers. Numerical tests show that each family includes schemes that are highly accurate compared to other mass-preserving methods that can be found in the literature.\\

}


\end{abstract}

\section{Introduction}\label{intro}

Conservation laws are among the most fundamental properties of a given system of partial differential equations (PDEs); all solutions must satisfy every conservation law. Consequently, there is a need to understand how to construct finite difference schemes that preserve multiple conservation laws of a given system.

For simplicity, we will limit our discussion to scalar PDEs for $u(x,t)$, where $x$ is a spatial variable and $t$ denotes time. (The generalization to systems with two or more independent variables is straightforward, but the notation becomes messier.) A local conservation law of a given PDE,
\[
\mathcal{A}(x,t,[u])=0,
\]
is a divergence expression,
\[
\text{Div}\,\mathbf{F}=D_x\{F(x,t,[u])\}+D_t\{G(x,t,[u])\},
\]
that is zero on the set of all solutions of the PDE. Here and henceforth, square brackets around a differentiable expression denote the expression and a finite number of its differential consequences using the total derivatives $D_x$ and $D_t$. In this notation,
\[\text{Div}\,\mathbf{F}=0 \quad \text{when}\quad [\mathcal{A}=0].\]
The functions $F$ and $G$ are, respectively, the flux and density of the conservation law. 

For suitable boundary conditions, the existence of such a conservation law implies that $\int G\upd x$ is an invariant, which means that it is constant on any solution. Some classes of PDEs have conservation laws and/or invariants built in as part of their structure. In particular, Hamiltonian PDEs are of the form
\begin{equation}\label{Hamstruct}
u_t=\mathcal{D}\frac{\delta}{\delta u}\mathcal{H},\qquad \mathcal{H}=\int H([u])\upd x,
\end{equation}
where $\mathcal{D}$ is a skew-adjoint differential operator, $\delta/\delta u$ denotes the variational derivative and $H$ is a given function. For every Hamiltonian PDE, the Hamiltonian functional $\mathcal{H}$ is an invariant.

Any Hamiltonian PDE can be embedded in a multisymplectic system; all such systems have a conservation law that expresses the invariance of a multisymplectic form \cite{Bridges,BHL}. Multisymplectic finite difference schemes, which are extensions of symplectic methods for Hamiltonian ordinary differential equations (ODEs), preserve this conservation law \cite{BridgesReich2001,BridgesReich,Leimku,Marsden,Ascher,AscherMcL,QuispMS}.

Sanz-Serna, de Frutos and Dur\'{a}n were among the first to study the benefits of using finite difference methods to preserve invariants of a Hamiltonian PDE \cite{DeFrutos,DuranLM,DuranSS1,DuranSS2}, using what is effectively a semi-discretization in time. An alternative approach is to discretize in space first, then apply an invariant-conserving method in time to preserve the approximated Hamiltonian functional. With this approach, invariants can be conserved by symplectic methods \cite{BridgesReich2001,Ascher,Cano,oliver}, discrete line integral methods \cite{BarlettiNLS,BI} or discrete gradient methods \cite{Celledoni,Dahlby,Furihata,Furibook}. 

Among the most useful discrete gradient methods for approximating Hamiltonian PDEs is the Average Vector Field (AVF) method introduced in \cite{Celledoni}, which generalizes Quispel and McLaren's energy-preserving approach for Hamiltonian ODEs \cite{QMcLaren}. According to \cite{Celledoni}, the AVF method ``is distinguished by its features of linear covariance, automatic preservation of linear symmetries, reversibility with respect to linear reversing symmetries and often by its simplicity.'' Here is an outline of this method for a given PDE of the form (\ref{Hamstruct}).

Relative to a generic lattice point $\mathbf{n}=(m,n)$, the uniform grid points are 
\begin{equation}\label{grid}
x_i=x(m+i)=x(m)+i\Delta x,\quad t_j=t(n+j)=t(n)+j\Delta t.
\end{equation}
In the following, $u$ is approximated by either the semidiscretization $U_i\approx  u(x_i,t)$ or the full discretization $u_{i,j}\approx u(x_i,t_j)$. The vectors whose $i$-th components are $U_i$ and $u_{i,j}$ are denoted by $\mathbf{U}$ and $\mathbf{u}_j$ respectively, and $\nabla$ is the differential operator whose $i$-th component is $\partial/\partial U_i$. 
Let $\widehat{\mathcal{H}}(\mathbf{U})\Delta x$ be a semidiscretization of $\mathcal{H}$ and let $\widetilde{\mathcal{D}}$ be a discrete skew-adjoint operator (typically depending on $\mathbf{u}_0$ and $\mathbf{u}_1$) that approximates $\mathcal{D}$.
The AVF method approximates (\ref{Hamstruct}) by 
\begin{equation}\label{AVF}
D_n \mathbf{u}_0=\widetilde{\mathcal{D}}\left(\int_0^1 \nabla \widehat{\mathcal H}\left(\xi \mathbf{u}_1+(1-\xi)\mathbf{u}_{0}\right)\mathrm{d}\xi\right).
\end{equation}

In \cite{McLQuisp}, McLachlan and Quispel proved that if the spatial discretization preserves a semidiscrete conservation law of the Hamiltonian PDE, the application of any discrete gradient method in time \cite{Gonz,Cell2,Robi,HLWbook,QuispTurn} yields a fully discrete (local) conservation law of the Hamiltonian. This is a powerful result, because conservation laws contain far more information than integral invariants do.

The current paper focuses on numerical methods that preserve two conservation laws of the modified Korteweg-de Vries (mKdV) equation,
\begin{equation}\label{mKdV}
\mca\equiv u_t+u^2u_x+u_{xxx}=0.
\end{equation}
This PDE possesses a bi-Hamiltonian structure: there are two different ways to write it in the Hamiltonian form \eqref{Hamstruct}, namely with
\begin{equation}\label{Hamiltonian1}
\mathcal{H}=-\int\tfrac{1}{12}u^4+\tfrac{1}2uu_{xx}\, \upd x,\qquad \mathcal{D}=D_x,
\end{equation}
and with
\begin{equation}\label{Hamiltonian2}
\mathcal{H}=-\int\tfrac{1}{2}u^2\,\mathrm{d}x,\qquad\mathcal{D}=D_x^3+\tfrac{1}{2}(u^2D_x+uu_x).
\end{equation}

The mKdV equation has infinitely many independent conservation laws \cite{Olverpaper}. Those for mass, momentum and energy, respectively, are \cite{Anco2,Miura}:
\begin{align}\label{mKdVCL1}
& D_xF_1+D_tG_1 = D_x\!\left(\tfrac{1}3 u^3 + u_{xx}\right)+D_tu,\\\label{mKdVCL2}
& D_xF_2+D_tG_2 =D_x\!\left(\tfrac{1}4u^4+uu_{xx}-\tfrac{1}{2}u_x^2\right)+D_t\!\left(\tfrac{1}2 u^2\right),\\\label{mKdVCL3}
& D_xF_3+D_tG_3 = D_x\!\left(\tfrac{1}2\!\left(\tfrac{1}3u^3+u_{xx}\right)^2+\tfrac{1}2u_xu_{t}-\tfrac{1}2uu_{xt}\right)+D_t\!\left(\tfrac{1}{12}u^4+\tfrac{1}2uu_{xx}\right).
\end{align}
For suitable (e.g. zero) boundary conditions, integrating (\ref{mKdVCL1})--(\ref{mKdVCL3}) over the spatial domain yields the invariants
\begin{equation}\label{glinvmkdv}
\int G_1\,\upd x=\int u\,\upd x,\quad \int G_2\,\upd x=\tfrac{1}{2}\int u^2\,\upd x, \quad \int G_3\,\upd x=\int \tfrac{1}{12}u^4+\tfrac{1}2uu_{xx}\, \upd x.
\end{equation}

For discretizations of the mKdV equation, the following notation is helpful. The forward shift operators in space and time are defined, for all functions $f$ on the grid, by
\[
S_m(f(x_i,t_j))=f(x_{i+1},t_j),\qquad S_n(f(x_i,t_j))=f(x_i,t_{j+1}).
\]
The forward difference operators $D_m, D_n$ and forward average operators $\mu_m,\mu_n$ are
$$D_m=\tfrac{1}{\Delta x}(S_m-I),\qquad D_n=\tfrac{1}{\Delta t}(S_n-I),\qquad\mu_m=\tfrac{1}{2}(S_m+I),\qquad \mu_n=\tfrac{1}{2}(S_n+I),$$
where $I$ is the identity operator.

The AVF method (\ref{AVF}) for (\ref{Hamiltonian1}) uses $$\widehat{\mathcal{H}}=-\sum_j\tfrac{1}{12}U_j^4+\tfrac{1}{2}U_jD_m^2U_{j-1},\qquad \widetilde{\mathcal{D}}=D_m\mu_mS_m^{-1}.$$
This gives the following 10-point energy-conserving scheme (and its shifts):
\[\text{AVF}_{\text{EC}}=D_m\mu_m\!\left\{\tfrac{1}{3}(\mu_nu_{-1,0})\,\mu_n\!\left(u_{-1,0}^2\right)+D_m^2\mu_nu_{-2,0}\right\}+D_nu_{0,0}.\]

To apply the AVF method (\ref{AVF}) to (\ref{Hamiltonian2}), let $\widehat{\mathcal{H}}=-\sum_j\tfrac{1}{2}U_j^2$ and let $\widetilde{\mathcal{D}}$ be the diagonal matrix whose entry $\widetilde{\mathcal{D}}_j$ in the row and column indexed by $j$ acts on functions $f(m,n,[u])$ as
\begin{equation*}
\begin{split}
\widetilde{\mathcal{D}}_jf=&\,D_m^3\mu_mS_m^{-2}f\!+\!\tfrac{1}{3}S_m^{-1}\mu_m\!\left\{\!\left((\mu_nu_{j,0})^2\!+\!(\mu_nu_{j,0})(\mu_nu_{j+1,0})\!+\!(\mu_nu_{j+1,0})^2\right)\!D_mf\right\}\\
&+\tfrac{1}{3}\left(\mu_nu_{j-1,0}+\mu_nu_{j,0}+\mu_nu_{j+1,0}\right)\left(D_m\mu_nu_{j-1,0}\right)f.
\end{split}
\end{equation*}
This yields the 10-point momentum-conserving scheme
\[\text{AVF}_{\text{MC}}=D_m\left\{\tfrac{1}{3}(\mu_m\mu_nu_{-1,0})\mu_m\big((\mu_nu_{-1,0})^2\big)+D_m^2\mu_m\mu_nu_{-2,0}\right\}+D_nu_{0,0}.\]

A new direct approach that enables multiple conservation laws to be preserved was introduced recently in \cite{Grant,GrantHydon} and greatly simplified in \cite{FCHydon}. This approach, which does not exploit either Hamiltonian structures or integrability, can be applied (at least, in principle) to any system of PDEs whose conservation laws are polynomial in $[u]$; see \S 4 of \cite{FCHydon} for a non-Hamiltonian, non-integrable example. In the next sections we show that the two AVF schemes $\text{AVF}_{\text{EC}}$ and $\text{AVF}_{\text{MC}}$ are particular members of two parametrized families of methods that can be found by using this approach.

The rest of this paper is organized as follows. Section~\ref{method} describes the simplified direct approach to finding finite difference schemes that preserve two local conservation laws, for mass and either momentum or energy. In \S\ref{mKdVsec}, we apply this strategy to the mKdV equation and obtain several families of conservative schemes. In \S\ref{test}, numerical tests confirm the theoretical results and demonstrate the effectiveness of the new schemes compared with multisymplectic and narrow box schemes, each of which preserves mass, but not momentum or energy.

\section{A strategy for preserving conservation laws}\label{method}

Given a (not necessarily Hamiltonian) PDE,
\begin{equation}\label{PDE}
\mathcal{A}(x,t,[u])=0,\qquad (x,t)\in V\subset\mathbb{R}^2,
\end{equation}
a conservation law is in characteristic form if it satisfies
\begin{equation}\label{charform}
\text{Div}\,\mathbf{F}=\mcq\mca.
\end{equation}
The function $\mcq$ is called the characteristic of the conservation law. For simplicity, we assume throughout that the domain $V$ is contractible.

\begin{rem}\label{Eulerrem} The vector space of total divergences is the kernel of the Euler operator, 
\begin{equation*}
\mathcal{E}=\sum_{i,j}(-D_x)^i(-D_t)^j\frac{\partial}{\partial u_{,ij}}\,,\qquad \text{where}\quad  u_{,ij}=D_x^iD_t^j(u).
\end{equation*}
\end{rem}
From Remark~\ref{Eulerrem}, if $$\mathcal{E}(\mcq\mca)=0,$$
there exists $\mathbf{F}$ such that $\text{Div}\,\mathbf{F}=\mcq\mca$ is a conservation law.

A finite difference approximation of (\ref{PDE}) on the grid is a partial difference equation,
\begin{equation}\label{PdE}
\wmca(m,n,[u_{0,0}]) =0,
\end{equation}
where square brackets $[\ ]$ around an expression defined on the grid denote the expression and a finite number of its shifts. From here on, tildes denote approximations to the corresponding continuous quantities. We will abbreviate $u_{0,0}$ to $u$ wherever this does not cause confusion.

The aim is to find finite difference schemes that satisfy a discrete analogue of each preserved conservation law,
\begin{equation}\label{discCLaw}
\mbox{Div}\,\widetilde{\mathbf{F}}=D_m\left\{\widetilde{F}(m,n,[u])\right\}+D_n\left\{\widetilde{G}(m,n,[u])\right\}=0,\quad\text{when}\quad [\wmca=0].
\end{equation}
This difference conservation law is in characteristic form if 
$$\mbox{Div}\,\widetilde{\mathbf{F}}=\wmcq(m,n,[u])\wmca(m,n,[u]);$$
just as for the continuous case, the multiplier function $\wmcq$ is called the characteristic (see \cite{Hydonbook} for more details).
The following result, due to Kuperschmidt \cite{Kuper} and generalized in \cite{HydonMans}, is a discrete version of Remark~\ref{Eulerrem} and plays a pivotal role in our approach.
\begin{rem}\label{discEulrem}
The kernel of the difference Euler operator
\[\mathsf{E}=\sum_{i,j}S_m^{-i}S_n^{-j}\frac{\partial}{\partial u_{i,j}},\]
is the vector space of all difference divergences (\ref{discCLaw}).
\end{rem}

Given a PDE (\ref{PDE}) with $p$ conservation laws in characteristic form  (\ref{charform}) that one wishes to preserve, our strategy is to seek approximations $\wmcq_i$ and $\wmca$ such that 
\begin{equation}\label{EQA}
\mathsf{E}(\wmcq_i\wmca)= 0,\qquad i=1,\dots p.
\end{equation}
From Remark~\ref{discEulrem}, there exist $\widetilde{\mathbf{F}}_i$ such that each $\wmcq_i\wmca=\text{Div}\,\widetilde{\mathbf{F}}_i$ is a discretization of the corresponding continuous conservation law. This strategy can be implemented efficiently as follows.
\begin{enumerate}
\item Choose a rectangular stencil that is large enough to contain second-order approximations of $\mca$ and all $\mcq_i$.
\item On the given stencil, the most general finite difference approximations $\wmca$ and $\wmcq_1$ will depend on a large number of coefficients.
\item Impose consistency conditions giving second-order accuracy of the approximations at the centre of the stencil (which need not be a grid point).
\item
Reduce the number of free parameters remaining by making some key terms as compact as possible; typically, these include highest derivatives and highest-order nonlinear terms.    
\item Some of the remaining parameters are determined by solving
\begin{equation}\label{Eulcond}
\mathsf{E}(\wmcq_1\wmca)=0
\end{equation}
symbolically, whenever a solution exists.
The discrete flux $\widetilde{F}_1$ and density $\widetilde{G}_1$, which satisfy $\wmcq_1\wmca=D_m \widetilde{F}_1+D_n\widetilde{G}_1$, can be reconstructed from the characteristics \cite{Hydon}.
\item Steps 2 onward are iterated to preserve more conservation laws. If $\mathsf{E}(\wmcq_i\wmca)\!=\!0$ has no solutions for some $i$, the corresponding conservation law cannot be preserved on the chosen stencil without violating at least one of the earlier conservation laws.
\end{enumerate}
Solving (\ref{EQA}) is a crucial step in the procedure; it is made tractable by the compactness conditions and the restriction to second-order approximations. Without these constraints, the only practical approach is to seek a Groebner basis (see \cite{Cox}) for the polynomial system that determines the coefficients. This can take weeks of symbolic computation, even for scalar equations whose conservation laws have nothing worse than quadratic nonlinearities, approximated on the smallest possible stencil. The complexity increases exponentially with the degree of the polynomial nonlinearity and the size of the stencil.

By contrast, the compactness and second-order conditions simplify the calculations to the point that schemes preserving multiple conservation laws can be found in just a few minutes; a Groebner basis may not even be needed. These simplifications were introduced in \cite{FCHydon} and used to obtain several parametrized families of schemes, some of which are highly accurate, for the KdV equation and a nonlinear wave equation with quadratic nonlinearity.

Recently, Frasca-Caccia  \cite{FC18} outlined a new one-parameter family of mass- and energy-conserving 10-point schemes for the mKdV equation, which has cubic nonlinearity. In the current paper, the simplified strategy is used to extend this result to schemes that preserve mass and either momentum or energy, using 8-point and 10-point stencils. This demonstrates that the simplified strategy is not limited to quadratic nonlinearity. 

\section{Conservative methods for the mKdV equation}\label{mKdVsec}

Each of the conservation laws (\ref{mKdVCL1})--(\ref{mKdVCL3}) is in characteristic form; their characteristics are
\begin{equation}\label{mKdVchar}
\mcq_1=1,\qquad \mcq_2=u,\qquad \mcq_3=\tfrac{1}3u^3+u_{xx}.
\end{equation}

\begin{figure}
\center{
\begin{tikzpicture}
\draw[fill] (0,3) circle [radius=0.075];
\draw[fill] (0,1.5) circle [radius=0.075];
\draw[fill] (2,1.5) circle [radius=0.075];
\draw[fill] (4,1.5) circle [radius=0.075];
\draw[fill] (6,1.5) circle [radius=0.075];
\draw[fill] (8,1.5) circle [radius=0.075];
\draw[fill] (2,3) circle [radius=0.075];
\draw[fill] (4,3) circle [radius=0.075];
\draw[fill] (6,3) circle [radius=0.075];
\draw[fill] (8,3) circle [radius=0.075];
\draw[thick,dotted] (0,1.5)--(8,1.5);
\draw[thick,dotted] (0,3)--(8,3);
\node [below right] at (3.9,1.45) {\scriptsize{$(x,t\!-\!{\Delta t/2})$}};
\node [below right] at (3.9,2.2) {\scriptsize{$(x,t)$}};
\node [below] at (3,2.2) {\scriptsize{$(x\!-\!{\Delta x/2},t)$}};
\draw (4,2.25) circle [radius=0.125];
\draw (3.85,1.35)--(4.15,1.65);
\draw (3.85,1.65)--(4.15,1.35);
\draw (2.85,2.10)--(3.15,2.4);
\draw (2.85,2.4)--(3.15,2.1);
\node [below] at (0,1.5) {\scriptsize{$A$}};
\node [below] at (8,1.5) {\scriptsize{$B$}};
\node [left] at (0,1.5) {\scriptsize{$C=0$}};
\node [left] at (0,3) {\scriptsize{$D=1$}};
\end{tikzpicture}
\caption{Example of a rectangular stencil for mKdV. PDEs and conservation laws are preserved to second order at the central point $(x, t)$; densities and fluxes are second-order at $(x, t-\Delta t/2) $ and $(x-\Delta x/2, t)$, respectively.}
\label{stencil}
}
\end{figure}
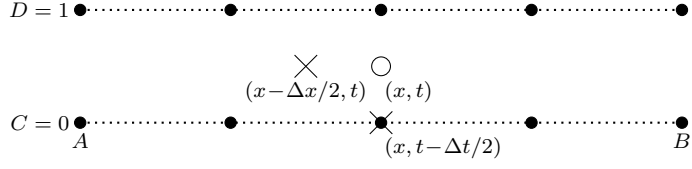

For simplicity, we consider only one-step schemes for the mKdV equation (\ref{mKdV}). Therefore, our stencils are as shown in Figure~\ref{stencil}, with $B\!-\!A\geq 4$.
On such stencils, we construct schemes for the mKdV equation of the form
\begin{equation}\label{standform}
\widetilde{\mca}=D_m\widetilde{F}_1+D_n\widetilde{G}_1=0,
\end{equation}
so that the mass conservation law (\ref{mKdVCL1}) is preserved.

From this starting-point, the strategy described in \S\ref{method} is used to preserve either (\ref{mKdVCL2}) or (\ref{mKdVCL3}). Linear terms in $F_1$, $G_1$ and $\mcq_2$ or $\mcq_3$ are approximated by linear combinations of the values of $u$ at the stencil points, with undetermined coefficients:
\begin{equation}\label{linterm}
\frac{\partial^{r+s}}{\partial x^r\partial t^s}u \approx \frac{1}{\Delta x^r}\frac{1}{\Delta t^s}\sum_{i=A}^B\sum_{j=0}^1\alpha_{i,j}u_{i,j}.
\end{equation}
Similarly, the cubic terms in $\mcq_3$ and $F_1$ are approximated by cubics:
\begin{equation}\label{cubterm}
\begin{split}
u^3\approx\sum_{i=A}^B&\left(\sum_{k=i}^B\sum_{j=k}^B(\beta_{i,j,k}u_{i,0}u_{j,0}u_{k,0}+\gamma_{i,j,k}u_{i,1}u_{j,1}u_{k,1})\right.\\
&\,\,\,\,\left.+\sum_{k=A}^B\sum_{j=k}^B(\delta_{i,j,k}u_{i,1}u_{j,0}u_{k,0}+\varepsilon_{i,j,k}u_{i,0}u_{j,1}u_{k,1})\right).
\end{split}
\end{equation}
The undetermined coefficients satisfy consistency conditions to ensure that the schemes are second-order accurate at the centre of the stencil.

\begin{rem}\label{rem2ord}
An approximation of a conservation law in the form (\ref{discCLaw}) is second-order accurate at the centre $(x, t)$ if the approximations of $\widetilde{F}$ and $\widetilde{G}$ are second-order accurate at the points $(x-\Delta x/2, t)$ and $(x, t-\Delta t/2)$, respectively.
\end{rem}

We now use the simplified strategy to obtain families of schemes, all of which depend on parameters that are $\mathcal{O}(\Delta x^2,\Delta t^2)$. It is possible to find values of the parameters that partially remove the leading terms of the local truncation error. However, there is no choice of the parameters that eliminates the second-order error terms identically; the optimal parameter values depend on the particular problem being approximated. No scheme in these families preserves all three conservation laws, though some come fairly close to doing so.
 
\subsection*{8-point schemes}
The most compact stencil for the mKdV equation has eight points. We choose $A=-2, B=1$ in Figure~\ref{stencil}, and seek second-order approximations of characteristics, densities and fluxes at $(-1/2,1/2)$, $(-1/2,0)$ and $(-1,1/2)$, respectively.
\subsubsection*{Energy-conserving methods.}
To simplify the symbolic calculation giving energy-conserving schemes on the 8-point stencil, we use the following approximations to the second derivatives in $F_1$ and $Q_3$, respectively:
$$u_{xx}\approx D_m^2\mu_nu_{-2,0}+\lambda D_mD_n\mu_m u_{-2,0} ,\,\, u_{xx}\approx D_m^2\mu_m\mu_nu_{-2,0}+\lambda D_mD_n\mu_m^2u_{-2,0},$$ where $\lambda=\mathcal{O}(\Delta x^2,\Delta t^2)$. All other terms in $\widetilde{\mca}$ and $\widetilde{Q}_3$ are of the form (\ref{linterm}) or (\ref{cubterm}), subject only to second-order consistency and
\begin{equation*}
\mathbf{E}(\widetilde{\mcq}_3\widetilde{\mca})= 0.
\end{equation*}
In this way, we find a family of schemes that depends on two parameters. One of these parameters is $\mathcal{O}(\Delta x^4,\Delta t^4)$; as this is negligible on fine grids, we set this parameter to zero, obtaining a family of schemes $\wmca(\lambda)$ with
\[
\widetilde{F}_1\!=\!\tfrac{1}3\!\left(\mu_n\mu_m^2u_{-2,0}\right)\!\mu_n\!\left(\!\left(\mu_m^2u_{-2,0}\right)^2\right)+D_m^2\mu_n u_{-2,0}+\lambda D_nD_m\mu_mu_{-2,0},\, \widetilde{G}_1\!=\!\mu_mu_{-1,0}.
\]
Each of these schemes preserves the following discrete version of the conservation law (\ref{mKdVCL3}):
\[
\widetilde{\mcq}_3\widetilde{\mca}= D_m\widetilde{F}_3+ D_n\widetilde{G}_3,
\]
with $\widetilde{\mcq}_3=\mu_m\widetilde{F}_1$ and
\begin{align*}\nonumber
\widetilde{F}_3=&\,\tfrac{1}2\widetilde{F}_1^2+\tfrac{1}2\left\{(D_m\mu_m\mu_n u_{-2,0}) D_n\mu_m^2 u_{-2,0}-(\mu_m^2\mu_nu_{-2,0})D_mD_n\mu_mu_{-2,0}\right\}\\
&+\tfrac{1}{2}\lambda(D_n\mu_mu_{-2,0})D_n\mu_mu_{-1,0}+\tfrac{1}{48}\Delta x^2\!\left\{(\mu_m^2u_{-2,0})^2+(\mu_m^2u_{-2,1})^2+(\mu_m^2u_{-2,0})(\mu_m^2u_{-2,1})\right\}\\
&\times\left\{(D_n\mu_m^2u_{-2,0})D_m\mu_m\mu_nu_{-2,0}-(\mu_m^2\mu_nu_{-2,0})D_mD_n\mu_mu_{-2,0}\right\},\\
\widetilde{G}_3=&\,\tfrac{1}{12}\left(\mu_mu_{-1,0}\right)\left(\mu_m^3u_{-2,0}\right)\left\{\mu_m^2\!\left(\left(\mu_mu_{-2,0}\right)^2\right)+\tfrac{1}{4}\Delta x^2(D_m\mu_mu_{-2,0})D_m\mu_mu_{-1,0}\right\}\\
&+\tfrac{1}2(\mu_mu_{-1,0})D_m^2\mu_mu_{-2,0}\,.
\end{align*}
In both  $\widetilde{G}_3$ and $\widetilde{F}_3$ there are terms that do not have any counterpart in the continuous density and flux. These all vanish as the stepsizes tend to zero.

Assuming that $\Delta t<\Delta x$, the leading error terms are $\mathcal{O}(\Delta x^2)$. We use the notation
$$\text{EC}_8(\lambda_1)=\wmca(\lambda_1\Delta x^2);$$
the parameter $\lambda_1$ may be chosen to minimise the local truncation error.
\subsubsection*{Momentum-conserving methods.}
The momentum-conserving schemes introduced in this section are obtained by specifying that
$$\widetilde{u_{xx}}=D_m^2\mu_nu_{-2,0},$$
and choosing the following compact approximations of $G_1$ and $Q_2\,$:
\begin{equation}\label{G1Q2MC8}
\widetilde{G}_1=\mu_m u_{-1,0},\qquad 
\widetilde{Q}_2=\mu_m\mu_nu_{-1,0}\,.
\end{equation}
Then $F_1$ is determined from the general forms (\ref{linterm}) and (\ref{cubterm}) by requiring consistency and 
\begin{equation*}
\mathbf{E}(\widetilde{\mcq}_2\widetilde{\mca})\equiv 0.
\end{equation*}
This gives a three-parameter family of schemes that preserve mass and momentum.\footnote{Without the constraint on $\widetilde{G}_1$ in (\ref{G1Q2MC8}), there is another parameter. Removing (instead) the constraint on $\widetilde{Q}_2$ yields two other families of schemes.} Two of the parameters are $\mathcal{O}(\Delta x^4,\Delta t^4)$, so we set these to zero and obtain the one-parameter family $\mathcal{A}(\lambda)$ with $\widetilde{G}_1$ given in (\ref{G1Q2MC8}) and
\begin{equation*}
\begin{split}
\widetilde{F}_1\!=&\,\tfrac{1}{6}\!\left(\mu_n(u_{-2,0}\!+\!u_{0,0})\right)\!\left(\mu_n u_{-1,0}\right)^2\!\!+\!D_m^2\mu_nu_{-2,0}\!+\!\lambda\!\left\{2(\mu_nu_{-1,0})\mu_m\!\left((D_m\mu_nu_{-2,0})^2\right)\right.\\
&\left.\!\!+2\!\left(D_m^2\mu_nu_{-2,0}\right)\!\mu_m^2\!\left((\mu_nu_{-2,0})^2\right)\!-\!\Delta x\Delta t\left(D_nD_m\mu_mu_{-2,0}\right)\!\mu_n\!\left((D_m\mu_mu_{-2,0})^2\right)\right\}
\end{split}
\end{equation*}
where $\lambda=\mathcal{O}(\Delta x^2,\Delta t^2)$.
For any value of $\lambda$, the discrete momentum conservation law 
$$\widetilde{\mcq}_2\widetilde{\mca}= D_m\widetilde{F}_2+ D_n\widetilde{G}_2=0$$
is preserved, with $\widetilde{\mathcal{Q}}_2$ given in (\ref{G1Q2MC8}) and
\begin{align*}
\widetilde{F}_2\!=&\,\tfrac{1}{12}(\mu_nu_{-1,0})^2\{2(\mu_nu_{-2,0})(\mu_m\mu_nu_{-1,0})+(\mu_nu_{-1,0})(\mu_nu_{0,0})\}+(\mu_nu_{-1,0})D_m^2\mu_nu_{-2,0}\\
&-\tfrac{1}2(D_m\mu_nu_{-2,0})D_m\mu_nu_{-1,0}+\lambda(\mu_m\mu_nu_{-2,0})(\mu_m\mu_nu_{-1,0})\left\{(D_m\mu_nu_{-2,0})D_m\mu_nu_{-1,0}\right.\\
&\left.+(\mu_n(u_{-2,0}+u_{0,0}))D_m^2\mu_nu_{-2,0}\right\}+\tfrac{1}4\lambda\Delta x\Delta t\left\{(D_m\mu_m\mu_nu_{-2,0})D_n\mu_m^2u_{-2,0}\right.\\
&\left.-(\mu_m^2\mu_nu_{-2,0})D_mD_n\mu_mu_{-2,0}\right\}\left\{2\mu_n\!\left((D_m\mu_mu_{-2,0})^2\right)-(D_m\mu_mu_{-2,0})D_m\mu_mu_{-2,1}\right\},\\
\widetilde{G}_2\!=&\,\tfrac{1}{2}(\mu_mu_{-1,0})^2\!\!+\!\lambda\Delta t\Delta x(\mu_mu_{-1,0})(D_m^2\mu_mu_{-2,0})\!\left\{\!\tfrac{1}4(D_m\mu_mu_{-1,0})D_m\mu_mu_{-2,0}\!-\!(D_m\mu_m^2u_{-2,0})^2\right\}\!.
\end{align*}
All terms in $\widetilde{F}_1$, $\widetilde{F}_2$ and $\widetilde{G}_2$ that do not approximate any quantity in the continuous expressions vanish as the stepsizes tend to zero.
They are identically zero if $\lambda=0$. 

If $\Delta t\! <\! \Delta x$, the leading term in the local truncation error is $\mathcal{O}(\Delta x^2)$; by choosing $\lambda = \lambda_2 \Delta x^2$ optimally, one may be able to remove at least part of this error. In the numerical tests section, we use the notation
\[\mbox{MC}_8(\lambda_2)=\widetilde{\mca}(\lambda_2\Delta {x}^2).\]

\subsection*{10-point schemes}
We now develop schemes that preserve two conservation laws of the mKdV equation on the 10-point stencil with $A=-2, B=2$ in Figure~\ref{stencil}. Approximations of characteristics, densities and fluxes of the conservation laws are second-order accurate at $(0,1/2)$, $(0,0)$ and $(-1/2,1/2)$, respectively.

\subsubsection*{Energy-conserving methods.}
To develop mass and energy conserving schemes on the 10-point stencil, we approximate $\widetilde{G}_1$ and the cubic term of $\widetilde{\mathcal{Q}}_3$ on the most compact possible sub-stencils. This gives the one-parameter family of schemes found in \cite{FC18}, namely $\wmca(\lambda)$ with $\widetilde{F}_1=\mu_m\varphi_{-1,0}$, $\widetilde{G}_1=u_{0,0}$, where
\begin{align*}
\varphi_{-1,0}&=\tfrac{1}3(\mu_n u_{-1,0}^2)(\mu_n u_{-1,0})+D_m^2\mu_n u_{-2,0}+\lambda D_mD_n\mu_mu_{-2,0},
\end{align*}
where $\lambda=\mathcal{O}(\Delta x^2,\Delta t^2)$. The discrete local energy conservation law satisfied by each scheme in this family is
\begin{align*}
\widetilde{\mcq}_3\widetilde{\mca}=& D_m\widetilde{F}_3+ D_n\widetilde{G}_3=0,
\end{align*}
where $\widetilde{\mcq}_3=\varphi_{0,0}$ and
\begin{align*}
\widetilde{F}_3&=\tfrac{1}2\big\{\varphi_{-1,0}\varphi_{0,0}+(D_m\mu_nu_{-1,0})D_n\mu_mu_{-1,0}-(\mu_m\mu_nu_{-1,0})D_mD_nu_{-1,0}+\lambda(D_nu_{0,0})D_nu_{-1,0}\big\},\\
\widetilde{G}_3&=\,\tfrac{1}{12}u_{0,0}^4+\tfrac{1}2u_{0,0}D_m^2u_{-1,0}\,.
\end{align*}
This family of schemes is written as $$\text{EC}_{10}(\lambda_3)=\wmca(\lambda_3\Delta x^2).$$
Note that EC$_{10}(0)$ is the AVF scheme AVF$_\text{EC}$ that was derived in \S\ref{intro}.

\subsubsection*{Momentum-conserving methods.}
The complexity of the symbolic computation for solving
$$\mathbf{E}(\widetilde{\mcq}_2\widetilde\mca)\equiv 0,$$
is reduced by using the most compact second-order approximations of $\widetilde{G}_1$ and $\widetilde{\mcq}_2$. This gives a family depending on 27 parameters. For simplicity, we discuss here only members of this family having a compact approximation of the nonlinear term in $\widetilde{F}_1$, ignoring parameters that are $\mathcal{O}(\Delta x^4,\Delta t^4)$. This yields a one-parameter family $\widetilde{\mca}(\lambda)$ with
\[
\widetilde{F}_1=\tfrac{1}3(\mu_m\mu_nu_{-1,0})\,\mu_m\!\left(\left(\mu_n u_{-1,0}\right)^2\right)\!+\!D_m^2\mu_n\mu_mu_{-2,0}\!+\!\lambda D_mD_nu_{-1,0}\,,\,\,\, \widetilde{G}_1=u_{0,0},
\]
where $\lambda=\mathcal{O}(\Delta x^2,\Delta t^2)$. The discrete local conservation law for momentum is
\begin{align*}
\widetilde{\mcq}_2\widetilde{\mca}=& D_m\widetilde{F}_2+ D_n\widetilde{G}_2=0,
\end{align*}
where $\widetilde{\mathcal Q_2}=\,\mu_n\mu_mu_{-1,0}$ and
\begin{align*}
\begin{split}
\widetilde{F}_2=&\,\tfrac{1}{12}(\mu_nu_{-1,0})(\mu_nu_{0,0})\!\left\{\left(\mu_nu_{-1,0}\right)^2+\left(\mu_nu_{0,0}\right)^2+(\mu_nu_{-1,0})(\mu_nu_{0,0})\right\}\\
&+(\mu_m\mu_nu_{-1,0})D_m^2\mu_m\mu_nu_{-2,0}-\tfrac{1}4(D_m\mu_nu_{-1,0})D_m\mu_n(u_{-2,0}+u_{0,0})\\
&+\tfrac{1}2\,\lambda\left\{(\mu_m\mu_nu_{-1,0})D_mD_nu_{-1,0}-(D_m\mu_nu_{-1,0})D_n\mu_mu_{-1,0}\right\},\\
\widetilde{G}_2=&\,\tfrac{1}2u_{0,0}\left(u_{0,0}+\lambda D_m^2u_{-1,0}\right).
\end{split}
\end{align*}
Again, all quantities that do not approximate any term in the corresponding continuous conservation law are identically zero when $\lambda=0$. 
We denote this family of schemes by $$\mbox{MC}_{10}(\lambda_4)=\widetilde{\mca}(\lambda_4\Delta x ^2).$$
The simplest scheme, MC$_{10}(0)$, is AVF$_\text{MC}$.

\section{Numerical tests}\label{test}

In this section we consider the mKdV equation (\ref{mKdV}) on a domain $V=[a,b]\times[0,T]$, with periodic boundary conditions. We use some benchmark tests to show the effectiveness of the schemes developed in \S\ref{mKdVsec} compared with two well-known conservative methods (see \cite{Aydin,Ascher}). The narrow box scheme is obtained by applying a standard finite volume discretization, giving
\begin{align*}
D_n(\mu_mu_{-1,0})+D_m\!\left(\tfrac{1}{3}(\mu_nu_{-1,0})^3+D_m^2\mu_nu_{-2,0}\right)=0.
\end{align*}
The second method is multisymplectic and amounts to
\begin{align*}
D_n\!\left(\mu_m^3u_{-2,0}\right)+D_m\!\left(\tfrac{1}{3}\,\mu_m\!\left((\mu_m\mu_nu_{-2,0})^3\right)+D_m^2\mu_nu_{-2,0}\right)=0.
\end{align*}
This compact one-step scheme is a more efficient version of the popular Preissmann scheme \cite{Aydin,AscherMcL}.
Each of these schemes is in divergence form, preserving a discrete version of the mass conservation law (\ref{mKdVCL1}).

For every test in this section, the computational time is roughly the same for all of the numerical schemes. Therefore, the main difference between schemes is the solution error at the final time $t=T$, evaluated as
\begin{equation}\label{errsol}
\left.\frac{\|u-u_{exact}\|}{\|u_{exact}\|}\right\vert_{t=T}.
\end{equation}
On a grid with $M$ points in space and $N$ points in time, we evaluate the error in the conservation laws by measuring the error in the global invariants in (\ref{glinvmkdv}) as follows:
\begin{align}\label{errcl}
\text{Err}_\ell=&\,\Delta x\max_{j=1,\ldots,N}\left|\sum_{i=1}^M \left(\widetilde{G_\ell}(x_i,t_j)-\widetilde{G_\ell}(x_i,t_1)\right)\right|,\qquad \ell=1,2,3.
\end{align}
Whenever $\widetilde{G}_2$ or $\widetilde{G}_3$ is not defined for one of the schemes considered here, the corresponding error is instead evaluated as
\begin{align}
\text{Err}_2&=\Delta x\max_{j=1,\ldots,N}\left|\sum_{i=1}^M  \tfrac{1}2(v_{i,j}^2-v_{i,0}^2)\right|,\label{errorCLa}\\ 
\text{Err}_3&=\Delta x\max_{j=1,\ldots,N}\left|\sum_{i=1}^M \left( \tfrac{1}{12}(v_{i,j}^4-v_{i,0}^4)+\tfrac{1}2v_{i,j}D_m^2v_{i-1,j}-\tfrac{1}2v_{i,0}D_m^2v_{i-1,0}\right)\right|,\label{errorCLb}
\end{align}
where $v_{i,j}=u(a+i\Delta x,\, j\Delta t)$ for schemes defined on the 10-point stencil and
\[
v_{i,j}=\tfrac{1}{2}\big(u(a+(i-1)\Delta x,\, j\Delta t)+u(a+i\Delta x,\, j\Delta t)\big)
\] for 8-point schemes.

For each numerical test and family of schemes, we state the parameter value that minimizes the solution error (\ref{errsol}). None of the schemes preserve three conservation laws, so we state the parameter values $\lambda_i$ that optimize the error in the unpreserved invariant given by (\ref{errorCLa}) or (\ref{errorCLb}). We also include the results for the simplest schemes, obtained by setting each $\lambda_i$ to zero. As EC$_{10}(0)\equiv$ AVF$_\text{EC}$ and MC$_{10}(0)\equiv$ AVF$_\text{MC}$, this enables comparison with the AVF schemes.

The first benchmark problem is the interaction of two solitons, with the exact solution
\begin{equation*}
u(x,t)=\frac{2\sqrt{6}\,\kappa(\sqrt{c_1}\cosh{\xi_2}+\sqrt{c_2}\cosh{\xi_1})}{(\kappa^2-1)+\kappa^2\cosh{(\xi_1-\xi_2)}+\cosh{(\xi_1+\xi_2)}}\,;
\end{equation*}
here
\begin{equation*}
\kappa=\frac{\sqrt{c_1}+\sqrt{c_2}}{\sqrt{c_1}-\sqrt{c_2}}\,,\qquad
\xi_1=\sqrt{c_1}(x-c_1t+d_1),\qquad \xi_2=\sqrt{c_2}(x-c_2t+d_2).
\end{equation*}
We set
\[c_1=2.5,\qquad c_2=0.5,\qquad d_1=12,\qquad d_2=2.5,\]
and solve this problem on $\Omega=[-20,20]\times[0,10]$, using step lengths $\Delta x=0.1$ and $\Delta t=0.025$.

The schemes $\mbox{EC}_{8}(\lambda_1)$, $\mbox{MC}_8(\lambda_2)$, $\mbox{EC}_{10}(\lambda_3)$ and $\mbox{MC}_{10}(\lambda_4)$ give the minimal solution error when $\lambda_1=1$, $\lambda_2=-0.077$, $\lambda_3=0.04$ and $\lambda_4=0.19$. The values $\lambda_1=-0.05$, $\lambda_2=-0.073$, $\lambda_3=0.20$ and $\lambda_4\ll 0$ minimize the error in the invariant that is not preserved by the scheme. As the modulus of $\lambda_4$ is very large, the corresponding term in MC$_{10}(\lambda_4)$ is not merely a perturbation; it undermines the correct behaviour of the solution. So for this problem, minimizing the error in energy turns out to be a poor criterion for choosing the parameter.

\begin{figure}[htp]
\begin{center}
	{\includegraphics[width=16.5cm,height=19cm]{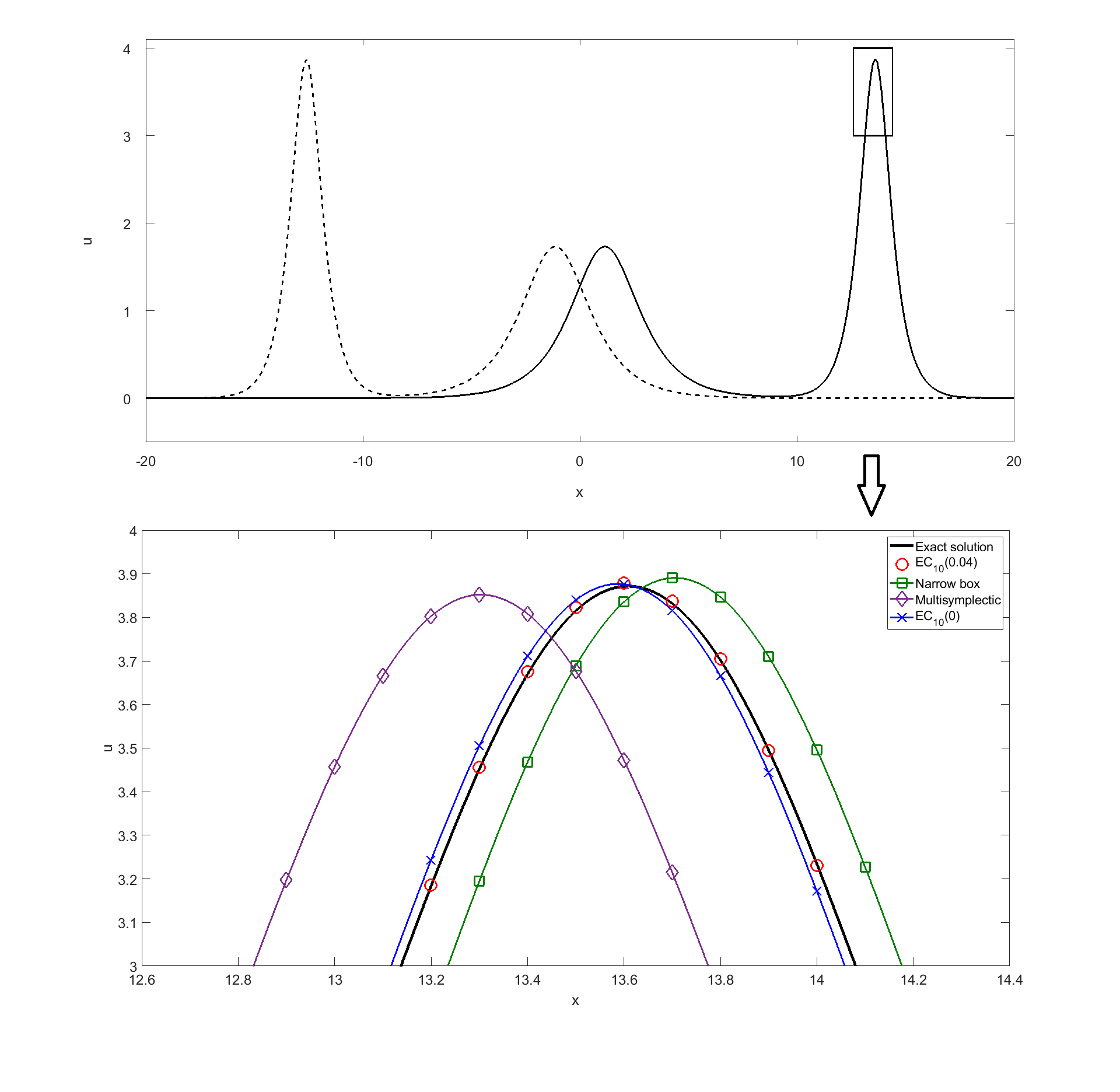}}
	\caption{Two-soliton problem for the mKdV equation. Top: 
		Initial condition (dashed line) and solution of $\mbox{EC}_{10}(0.04)$ with $\Delta x=0.1$, $\Delta t=0.025$ at time $T=10$ (solid line). Bottom: Top of the faster soliton; exact profile (solid line) and solutions of $\mbox{EC}_{10}(0.04)$ (circles), narrow box (squares), multisymplectic (diamonds) and $\mbox{EC}_{10}(0)$ (crosses).}
	\label{2solfig}
	\end{center}
\end{figure}

\begin{table}[t!]
	\caption{Errors in conservation laws and solutions for the two-soliton  problem for the mKdV equation, with $\Delta x\!=\!0.1$, $\Delta t\!=\!0.025$. An asterisk denotes the error that is minimized.}
	\label{2solfine}
	\small
	\centerline{\begin{tabular}{|c|c|c|c|c|c|c|c|c|}
			\hline
			\hline
			Method &  $\text{Err}_1$ & $\text{Err}_2$ & $\text{Err}_3$  & Sol. Err. & $\text{Err}_{\phi_1}$& $\text{Err}_{\phi_2}$&$\text{Err}_{\phi}$\\ 
			\hline
			\hline
			$\mbox{EC}_{8}(0)$	& 1.74e-13  & 0.0036 & 5.13e-13 &  0.3701 & -0.51 & -0.06 & -0.45\\
			\hline
			$\mbox{EC}_{8}(1)$	& 1.33e-13  & 0.0732 & 8.01e-13 & \phantom{$^*$}0.0085$^*$ & 0& -0.01 & 0.01 \\
			\hline
			$\mbox{EC}_{8}(-0.05)$	& 6.22e-14 & \phantom{$^*$}1.81e-04$^*$ & 4.65e-13 & 0.3857 & -0.53 & -0.07 & -0.46\\	
			\hline
			$\mbox{MC}_{8}(0)$	& 2.13e-13  & 3.69e-13 & 0.0632 &  0.2396 & -0.32 & -0.04 & -0.28\\
			\hline
			$\mbox{MC}_{8}(-0.077)$	& 1.21e-13  & 3.32e-13 & 0.0032 & \phantom{$^*$}0.0051$^*$& 0 & 0.01 & -0.01\\
			\hline
			$\mbox{MC}_{8}(-0.073)$	& 6.93e-14  & 1.46e-13 & \phantom{$^*$}5.55e-04$^*$ & 0.0139 & -0.02 & 0.01 & -0.03\\	
			\hline
			$\mbox{AVF}_\text{EC}$; $\mbox{EC}_{10}(0)$ & 3.91e-14  & 0.0142 & 4.80e-14 &  0.0167 & -0.02 & 0.01& -0.03\\
			\hline
			$\mbox{EC}_{10}(0.04)$	& 3.73e-14  & 0.0114 & 9.41e-14 & \phantom{$^*$}0.0030$^*$ & 0 & 0.01 & -0.01\\
			\hline
			$\mbox{EC}_{10}(0.20)$	& 5.15e-14 & \phantom{$^*$}1.82e-04$^*$ & 5.51e-14 & 0.0627 & 0.08& 0.02 &0.06\\
			\hline	
			$\mbox{AVF}_\text{MC}$; $\mbox{MC}_{10}(0)$	& 4.62e-14 & 5.68e-14 & 0.0358 & 0.0756 & -0.10 & 0 & -0.10\\
			\hline
			$\mbox{MC}_{10}(0.19)$	& 4.26e-14  & 5.33e-14 & 0.0359 & \phantom{$^*$}0.0051$^*$ & 0& 0.01 & -0.01\\
			\hline
			Narrow box & 1.28e-13  & 0.0117 & 7.0014 &  0.0742 & 0.10 & 0.02 & 0.08\\
			\hline
			Multisymplectic & 6.04e-14  & 0.0058 & 6.8991 &  0.2279 & -0.31 & -0.04 &-0.27\\
			\hline
			\hline
	\end{tabular}}
\end{table}

\begin{table}[b!]
	\caption{Errors in conservation laws and solutions for the two-soliton  problem for the mKdV equation, with $\Delta x=0.2$, $\Delta t=0.05$. An asterisk denotes the error that is minimized.}
	\label{2solcoarse}
	\small
	\centerline{\begin{tabular}{|c|c|c|c|c|c|c|c|c|}
			\hline
			\hline
			Method &  $\text{Err}_1$ & $\text{Err}_2$ & $\text{Err}_3$  & Sol. Err. & $\text{Err}_{\phi_1}$& $\text{Err}_{\phi_2}$& $\text{Err}_{\phi}$\\
			\hline
			\hline
			$\mbox{EC}_{8}(0)$	& 4.62e-14  & 0.0155 & 6.93e-14 &  0.9599 & -1.84 & -0.26 & -1.58\\
			\hline
			$\mbox{EC}_{8}(0.97)$	& 3.55e-14  & 0.2754 & 1.14e-13 & \phantom{$^*$}0.0358$^*$ & 0 & -0.03 & 0.03\\
			\hline
			$\mbox{EC}_{8}(-0.06)$	& 4.26e-14 & \phantom{$^*$}5.19e-04$^*$ & 1.15e-13 & 0.9798 & -1.93 & -0.26 & -1.67\\	
			\hline
			$\mbox{MC}_{8}(0)$	& 4.44e-14  & 9.41e-14 & 0.2363 &  0.7553 & -1.21 & -0.15 & -1.06\\
			\hline
			$\mbox{MC}_{8}(-0.079)$	& 6.57e-14  & 1.42e-13 & 0.0138 & \phantom{$^*$}0.0215$^*$& 0 & 0.05 & -0.05\\
			\hline
			$\mbox{MC}_{8}(-0.075)$	& 4.09e-14  & 7.11e-14 & \phantom{$^*$}0.0021$^*$ & 0.0567 & -0.06 & 0.04 & -0.10\\
			\hline
			$\mbox{AVF}_\text{EC}$; $\mbox{EC}_{10}(0)$ & 2.13e-14  & 0.0574 & 3.73e-14 &  0.0725 & -0.1 & 0.02 & -0.12\\	
			\hline
			$\mbox{EC}_{10}(0.05)$	& 2.66e-14  & 0.0438 & 4.09e-14 & \phantom{$^*$}0.0116$^*$ & 0 & 0.03 & -0.03\\
			\hline
			$\mbox{EC}_{10}(0.21)$	& 2.13e-14 & \phantom{$^*$}6.74e-04$^*$ & 4.97e-14 & 0.2571 & 0.35 & 0.08 & 0.27\\
			\hline	
			$\mbox{AVF}_\text{MC}$; $\mbox{MC}_{10}(0)$	& 1.95e-14 & 4.80e-14 & 0.1461 & 0.2959 & -0.40 & -0.01 & -0.39\\ 
			\hline
			$\mbox{MC}_{10}(0.19)$	& 2.31e-14  & 2.49e-14 & 0.1477 & \phantom{$^*$}0.0205$^*$ & 0 & 0.08 & -0.08\\
			\hline
			Narrow box & 4.97e-14  & 0.0459 & 6.8421 &  0.3054 & 0.40 & 0.09 & 0.31\\
			\hline
			Multisymplectic & 2.66e-14  & 0.0228 & 6.4635 &  0.7278 & -1.15 & -0.17 & -0.98\\
			\hline
			\hline
	\end{tabular}}
\end{table}

\clearpage

\begin{figure}[htp]
\begin{center}
	{\includegraphics[width=16.5cm,height=9cm]{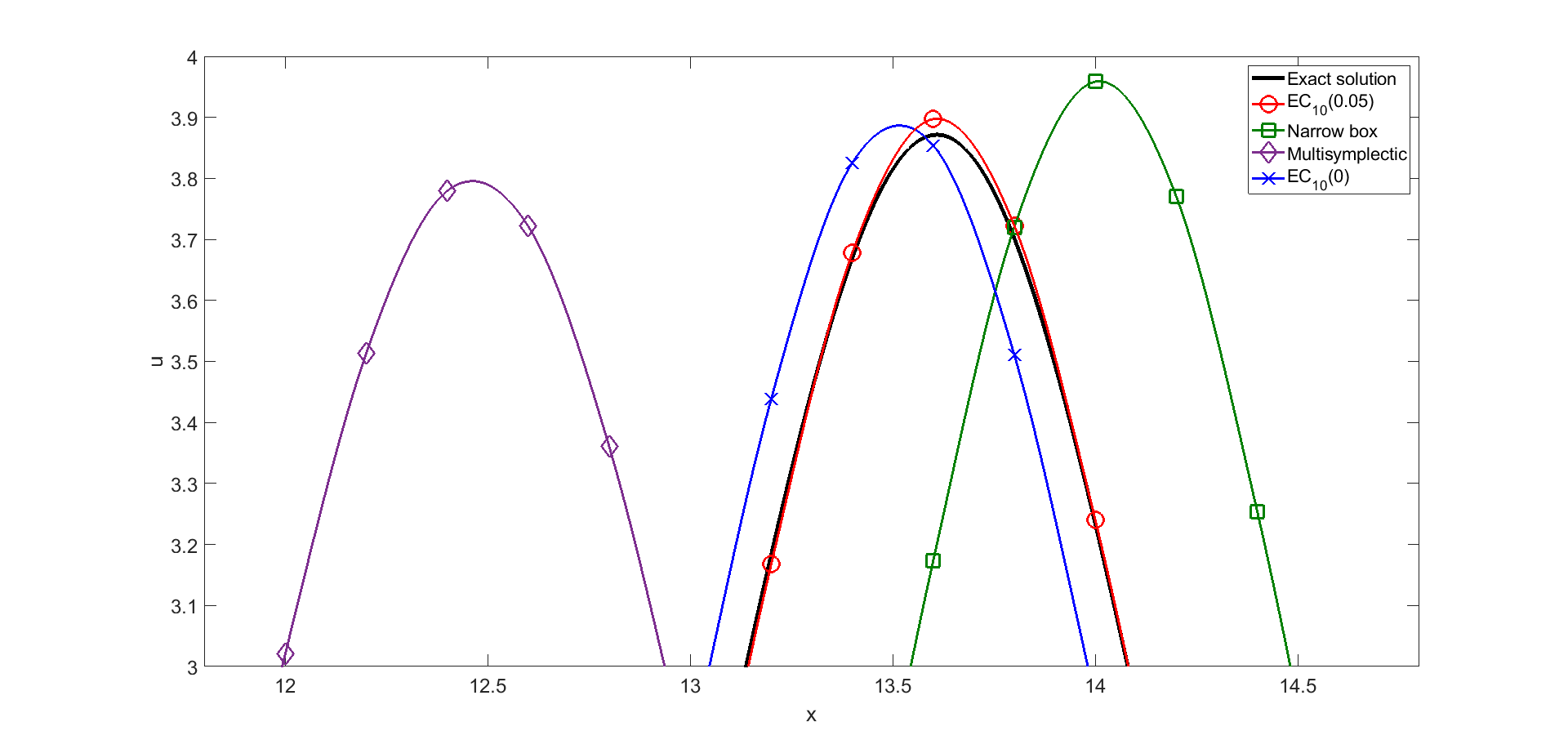}}
	\caption{Two-soliton problem for the mKdV equation  with $\Delta x=0.2$, $\Delta t=0.05$ at time $T=10$. Top of the soliton: exact profile (solid line) and solutions of $\mbox{EC}_{10}(0.05)$ (circles), narrow box scheme (squares), multisymplectic scheme (diamonds) and $\mbox{EC}_{10}(0)$ (crosses).}
	\label{2solfig2}
	\end{center}
\end{figure}

Table~\ref{2solfine} shows the errors at the final time $T=10$ in the conservation laws and the solution. We also show the error in the phase shift of the fastest and the slowest soliton, evaluated as, respectively,
$$\left.\text{Err}_{\phi_1}=(x_1-\tilde{x}_1)\right\vert_{t=10}, \qquad \left.\text{Err}_{\phi_2}=(x_2-\tilde{x}_2)\right\vert_{t=10};$$
where $x_{i}$ (resp. $\tilde{x}_{i}$) denotes the location of the soliton peak in the exact and (resp. numerical) solution, using piecewise cubic interpolation of approximations at the grid points.
The quantity 
\[\text{Err}_\phi= \text{Err}_{\phi_1}-\text{Err}_{\phi_2},\]
measures the extent to which the numerical solution underestimates the phase shift produced by the interaction of the two solitons.

Table~\ref{2solfine} shows that the new conservative schemes each preserve two discrete invariants (up to rounding errors). Each family includes schemes that are highly accurate --- more so than the multisymplectic and narrow box schemes. The small values of Err$_{\phi_1}$, Err$_{\phi_2}$ and Err$_{\phi}$ indicate that the best schemes also reproduce the correct phase shifts. Attempting to minimize the error in the unpreserved conservation law does not optimize the numerical solution; nevertheless, this approach gives MC$_8$ and EC$_{10}$ schemes that are more accurate than the narrow box and multisymplectic schemes. This is not true for EC$_8$ and MC$_{10}$ schemes.

The upper plot in Figure~\ref{2solfig} shows the initial condition and the numerical solution at $T=10$ given by the most accurate of our schemes, $\mbox{EC}_{10}(0.04)$. The lower plot compares various numerical solutions with the exact solution at $T=10$, near to the top of the faster soliton. The approximations at the grid points are connected by piecewise cubic interpolation, except for the solution of $\mbox{EC}_{10}(0.04)$ (for which the interpolation would cover the exact solution). The narrow box scheme is quite accurate for this problem, but it is outclassed by the AVF scheme $\mbox{EC}_{10}(0)$ and more so by $\mbox{EC}_{10}(0.04)$.

We now study the same problem on a coarser grid with $\Delta x=0.2$ and $\Delta t=0.05$. For each family, the parameter values minimizing the error in the solution or the unpreserved invariant are as given in Table~\ref{2solcoarse}; they are close to the optimal values for the finer grid. The solution error for the most accurate scheme in each family is around 4 times greater than on the finer grid, as expected. Figure~\ref{2solfig2} compares the exact solution with various numerical solutions near to the top of the fastest soliton, again using piecewise cubic interpolation.

As a second benchmark problem, we approximate the breather solution (see \cite{Zhang}),
\begin{equation*}
u(x,t)=\frac{\partial}{\partial x}
\left[-2\sqrt{6}\arctan\left(\frac{\sqrt{3}\sin{(2x-64t-\pi/2)}}{\cosh{(2\sqrt{3}x)}}\right)\right],
\end{equation*}
on the domain $\Omega=[-2,2]\times[0,0.4]$, setting $\Delta x=0.02$ and $\Delta t=0.002$.

Table~\ref{breathertable} shows the error in the conservation laws and the solution at the final time: $\mbox{EC}_{8}(2.22)$ is the most accurate scheme, while $\mbox{MC}_{8}(-0.165)$, $\mbox{EC}_{10}(0.92)$ and $\mbox{MC}_{10}(1.15)$ are the best in their respective families at minimizing the solution error. The error in the unpreserved conservation law is minimized by choosing $\lambda_1=0.49$, $\lambda_2=-0.128$, $\lambda_3=0.78$ and $\lambda_4\ll 0$. As for the two-soliton problem, minimizing the error in the non-conserved invariant is a poor criterion for choosing the parameter in MC$_{10}$, so we do not show this result.
By contrast, choosing the values of the parameters in $\mbox{MC}_8$ and $\mbox{EC}_{10}$ that minimize the error in the unpreserved invariant yields fairly accurate approximations.
\begin{figure}[htp]
\begin{center}
	{\includegraphics[width=16.5cm,height=9.5cm]{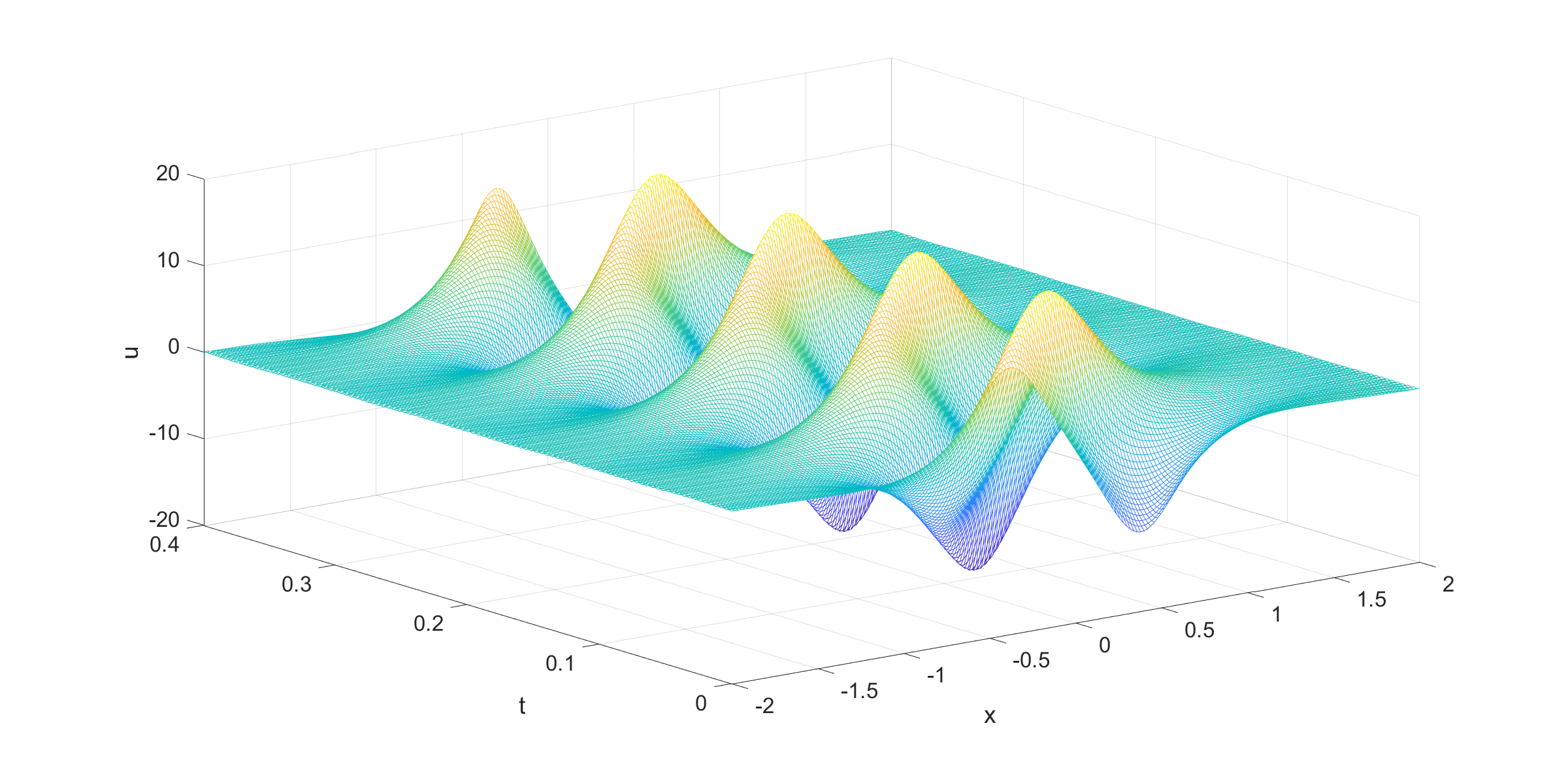}}
	{\includegraphics[width=16.5cm,height=9.5cm]{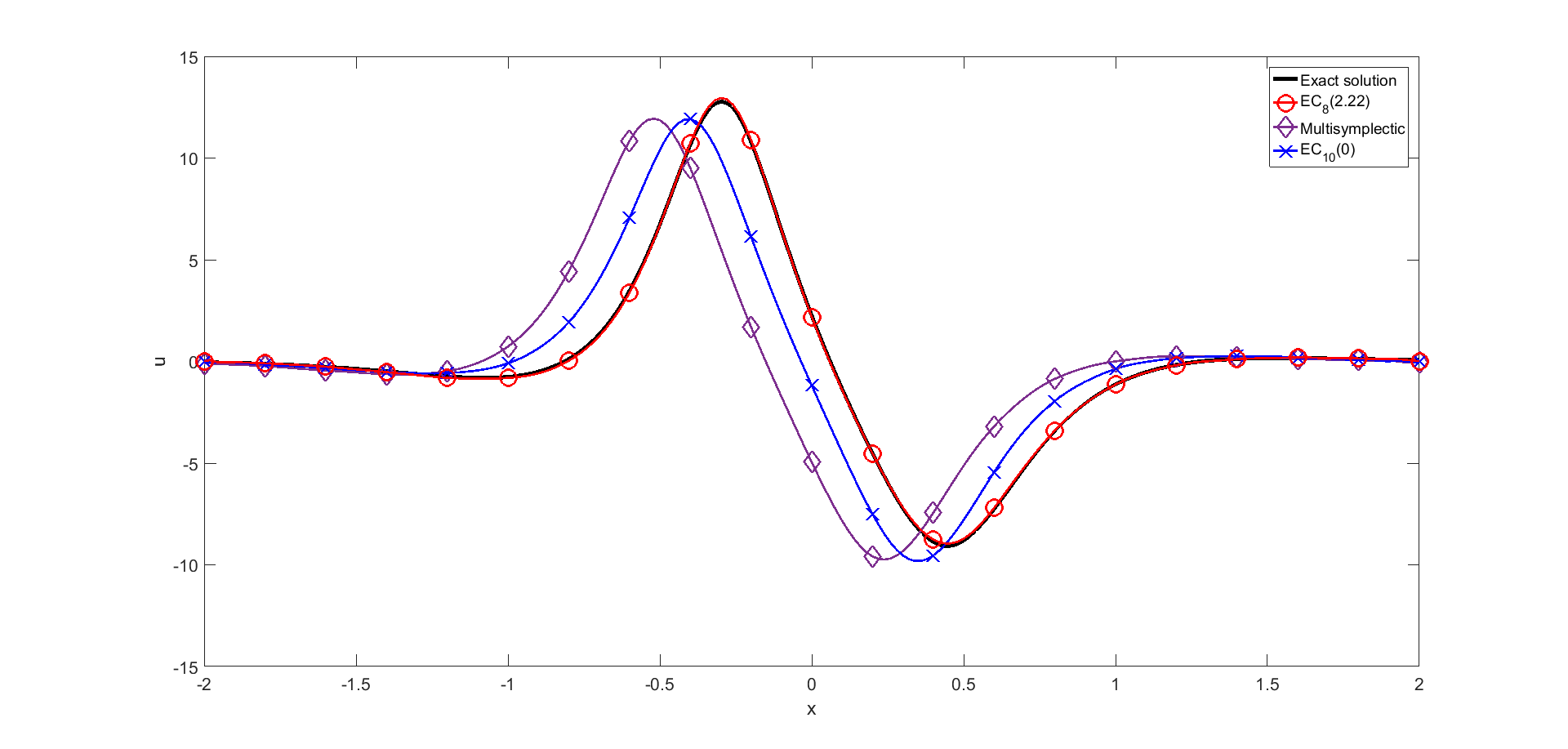}}
	\caption{Breather problem. Top: 
		Numerical solution given by method $\mbox{EC}_{8}(2.22)$ with stepsizes $\Delta x=0.02$, $\Delta t=0.002$ at time $T=0.4$. Bottom: exact profile (solid line) and solutions of methods $\mbox{EC}_{8}(2.22)$ (circles), multisymplectic (diamonds) and $\mbox{EC}_{10}(0)$ (crosses) at the final time (markers at every tenth point).}
	\label{breatherfig}
	\end{center}
\end{figure}

\begin{table}[t]
\caption{Errors in conservation laws and solution for the breather problem, setting $\Delta x=0.02$, $\Delta t=0.002$. An asterisk denotes the error that is minimized.}
\label{breathertable}
\small
\centerline{\begin{tabular}{|c|c|c|c|c|c|c|}
\hline
\hline
Method &  $\text{Err}_1$ & $\text{Err}_2$ & $\text{Err}_3$  &  Solution error \\ 
\hline
\hline
$\mbox{EC}_{8}(0)$	& 6.76e-13  & 0.1091 & 1.33e-10 &  0.9099\\
\hline
$\mbox{EC}_{8}(2.22)$	& 2.97e-13  & 0.3979 & 1.55e-10 & \phantom{$^*$}0.0144$^*$\\
\hline
$\mbox{EC}_{8}(0.49)$	& 9.59e-13 & \phantom{$^*$}0.0079$^*$ & 1.05e-10 & 0.7442\\	
\hline
$\mbox{MC}_{8}(0)$	& 6.25e-13  & 5.31e-12 & 7.534 & 0.7666  \\
\hline
$\mbox{MC}_{8}(-0.165)$	& 3.03e-13  & 2.74e-12 & 2.3599 & \phantom{$^*$}0.0497$^*$\\
\hline
$\mbox{MC}_{8}(-0.128)$	& 7.24e-13  & 7.60e-12 & \phantom{$^*$}0.1728$^*$ & 0.1931\\
\hline
$\mbox{AVF}_\text{EC}$; $\mbox{EC}_{10}(0)$ & 3.28e-14  & 0.1765 & 1.24e-11 &  0.4042\\
\hline
$\mbox{EC}_{10}(0.92)$	& 3.53e-14  & 0.0296 & 2.63e-11 & \phantom{$^*$}0.0295$^*$ \\
\hline
$\mbox{EC}_{10}(0.78)$	& 5.96e-14 & \phantom{$^*$}0.0095$^*$ & 1.35e-11 & 0.0708\\
\hline	
$\mbox{AVF}_\text{MC}$; $\mbox{MC}_{10}(0)$	& 1.14e-13 & 9.24e-13 & 4.3586 & 0.5040 \\
\hline
$\mbox{MC}_{10}(1.15)$	& 5.34e-14  & 7.03e-13 & 4.8298 & \phantom{$^*$}0.0219$^*$ \\
\hline
Narrow box & 1.00e-12  & 0.0382 & 566.37 &  0.3477 \\	
\hline
Multisymplectic & 2.60e-13  & 0.0184 & 539.40 &  0.7994\\
\hline
\hline
\end{tabular}}
\end{table}

The upper part of Figure~\ref{breatherfig} shows the numerical solution given by $\mbox{EC}_{8}(2.22)$. The lower part compares the exact solution and the numerical solutions given by EC$_{8}$(2.22), the multisymplectic scheme and EC$_{10}$(0) (which is the most accurate AVF scheme for this problem) at the final time. The graph of the solution of the narrow box scheme is very close to the solution of EC$_{10}(0)$, so is not shown in the figure. These schemes are more accurate than the multisymplectic scheme, but the frequency of the breather oscillations is best caught by EC$_{8}(2.22)$, whose solution almost overlaps the exact profile.
\section{Conclusions}\label{Conclusions}

The approach introduced in \cite{FCHydon}, which uses a fast symbolic computation to find finite difference schemes that preserve two conservation laws, is not restricted to quadratic nonlinearity. By considering stencils with eight and ten nodes, we have introduced four new one-parameter families of schemes that preserve two conservation laws of the mKdV equation, which has a cubic nonlinearity.

The AVF schemes introduced by Quispel and co-workers are found by this approach; typically, each is the simplest member of its family. However, the need to obtain a skew-adjoint approximation of the skew-adjoint differential operator means that AVF methods are restricted to stencils with an odd number of points in space. 

Each of the families includes schemes that are very accurate. However, none of them preserve the first three conservation laws. Nevertheless, the value of the parameter can be chosen to minimize the error in the third invariant. Our numerical tests have shown that this criterion leads to reasonably accurate solutions for some, but not all, families.

\section*{Acknowledgments} The authors would like to thank the Isaac Newton Institute for Mathematical Sciences for support and hospitality during the programme 
\textit{Geometry, Compatibility and Structure Preservation in Computational Differential Equations},
when work on this paper was undertaken. This work was supported by EPSRC grant number EP/R014604/1.

\end{document}